**STRUCTURE OF THE PRIME DISTRIBUTION WITHIN THE SEQUENCE OF NATURAL NUMBERS**

Andrei V. Vityazev and Galina V. Pechernikova*

**Abstract**

In the course of studies of the measure of chaos for the distribution of the prime numbers among the positive integers *N* arched structures have been found. It is given a brief description of the fine structure of the positive integers sequence, including the distribution of the prime numbers. A new non-Eratostenian sieve was built to demonstrate the multi-periodicity in the structure of the positive integers sequence. This sieve presents an infinite triangular matrix, which is named by us as R-matrix (Russian matrix). The simple deterministic equations for prime numbers and the main term of asymptotic were derived with the use of this matrix.

## Introduction

An interest in classical examples of irregular numeric sequences increases now in connection with the problem of organizations and chaos of quantum systems (in the broad sense) [1−5]. In the course of studies of the measure of chaos for the distribution of the prime numbers among the positive integers *N* we found [6, 7] arched structures. Here we discuss their nature and, in fact, we give a brief description of the fine structure of the positive integers sequence, including the distribution of the prime numbers. A new non-Eratostenian sieve with the memory is built to demonstrate the multi-periodicity in the structure of the positive integers sequence. This sieve presents an infinite triangular matrix, which is named by us as R-matrix (Russian matrix) [8] reflecting simultaneously the value and the unique "coloring" of each integer. With the use of the sieve we derived simple deterministic equations for prime numbers and the main term of asymptotic. Our results, in particular, mean that many natural apparent chaotic dependencies and space–time sequences appear wholly regular in coordinates of the system of prime numbers.

Notions of the positive integers sequence $N = \{1, 2, ...\}$ and system of natural numbers $<N, +, \times, 1>$ play a fundamental role in operations of ordering, measurement and calculation. The integers $p > 1$ having only two divisors 1 and $p$ form a peculiar frame of the positive integers sequence. These numbers are named prime numbers, but the other numbers in *N*, presented uniquely as a product of primes are named composite. The gaps between successive primes at the average increase with the growing $p$ but the combination of some regularity with elements of the chaos in distribution $p$ within *N* stays unclear.

The study of the distribution of primes among the positive integers with the problem of the solution for diophantine equations occupies a central place in Gilbert's 8-th problem. For physicists the prime numbers begin to play a fundamental role in base models of quantum systems, when the notion of the indivisible element is introduced. Vagueness of the distribution of primes among *N* and in other more complex numeric sequences, delivering headache to mathematicians and cryptologists, worries physicists as well. In the last few years one hypothesis is broadly discussed in connection with the problem of the quantum chaos. This hypothesis concerns possible analogies between the distribution of primes in numeric sequences and particularities of spectra of energy levels of excited quantum systems (atoms and nuclei in strong magnetic fields) [1−5]. Moreover, a searching for quantum-mechanical system is going in hope to check the Riemann hypothesis for all nontrivial zeroes of ζ-functions in complex planes are upon Re(*s*)=1/2. Situation aggravated even more in connection with the beginning of development of quantum computer technologies. It explains perpetual attempts to resolve the problem of the distribution of primes in *N*.

* Institute for Dynamics of Geospheres, Russian Academy of Sciences, 38 Leninsky prosp. (bldg. 1), 119334 Moscow, Russia, email: avit@idg.chph.ras.ru and pechernikova@idg.chph.ras.ru

Our approach to the problem differs in its strategy from other researches. Taking into account classical results for the first members of asymptotics of the distribution of primes, we considered behavior of the deviations from asymptotic and found arched structures unknown earlier [6, 7]. Below we give their brief description within the framework of elementary notions and theorems of the classical number theory and introduce an algebraic construction of nonEratosthenian sieve (R-matrix) [6-8].

**Necessary determinations and classical results** [9–12] are listed below:

- If $x$ is a positive real number, $\pi(x)$ is the number of primes less than or equal to $x$. As for example if primes under 25 are 2, 3, 5, 7, 11, 13, 17, 19, 23 then $\pi(3) = 2$, $\pi(10) = 4$, $\pi(25) = 9$.
- According to the theorem on the primes the main asymptotic is $\pi(x) \sim x / \ln x$, early experimentally found by Legendre in 1798.
- $p(n)$ is the $n$-th prime number.
- From asymptotic for $\pi(x)$ one can obtain $p(n) \sim n \ln n$.
- An improved evaluation gives $p(n) \sim n (\ln n + \ln \ln n - 1)$.
- Gauss' offer is $\pi(x) \sim \mathrm{li}(x)$.
- The consequence from the Hadamard's and de Valle-Poussen's theorem on primes is $\pi(x) = \mathrm{li}(x) + O(x \exp(- a (\ln x)^{1/2}))$.
- According to Koch the acceptance of Riemann' hypothesis is equivalent by $\pi(x) = \mathrm{li}(x) + O(x^{1/2} \ln x)$.
- $\Delta(n) = p(n) - \text{asymptotic}[p(n)]$

**I**. Let us consider the distributions of $\Delta(n)$ (Fig. 1, 2).

Figure 1 shows the distribution of the deviations of the primes $p(n)$ from their asymptotic values $n(\ln n + \ln \ln n - 1)$. Panel **a** displays the distribution of the spots $\Delta(n)$ in the range $n$ from 1 to 700. Initial area of this distribution under $n = 150$ is shown in **b**, where arched structures in the field of $50<n<140$ are seen distinctly. They can be seen and on the panel **a** as well. The distribution of $\Delta(n)$ for values $n$ in the interval 2200–3200 is shown on panel **c**, where our eye, accustomed to arched structures, distinguishes them everywhere.

Figure 2 shows the fluctuations $\Delta_1(n)$. The panel **a** displays a spline of the distribution $\Delta_1(n)$ obtained by the subtraction from $p(n)$ of its least squares numerical approximation on this interval. Fluctuations about the mean (about the abscissa axis $n$) grow at the average $\sim n^{1/2}$. The initial area the same curves with spots, forming up in arcs, is shown in **b** in increased scale.

First only the separate arcs can be stated, but a more careful consideration makes possible to state, that all the prime numbers belong to arcs, forming arched structures. Is it possible to reveal the arcs by examining spectrum $p(n)$? The Fourie-Walsh spectra of flattened distribution $p(n)$ and their asymptotic for sufficiently large $n$ were built. The conclusion is that having only the spectrum of $p(n)$, to catch arched structure it would be hardly possible. Probably, it explains that the references on arched structures in numeric sequences are full absent in the literature. The "best" approximation of $p(n)$ was used in order to exclude of the influence of asymptotic. Corresponding splines turn out to be closely symmetrical: the dispersion grows with the growing $n$ but asymmetry has not distinct trend with the growing $n$, as has been shown for the beginning of sequence on the Fig. 2a. What is the nature of these arches?

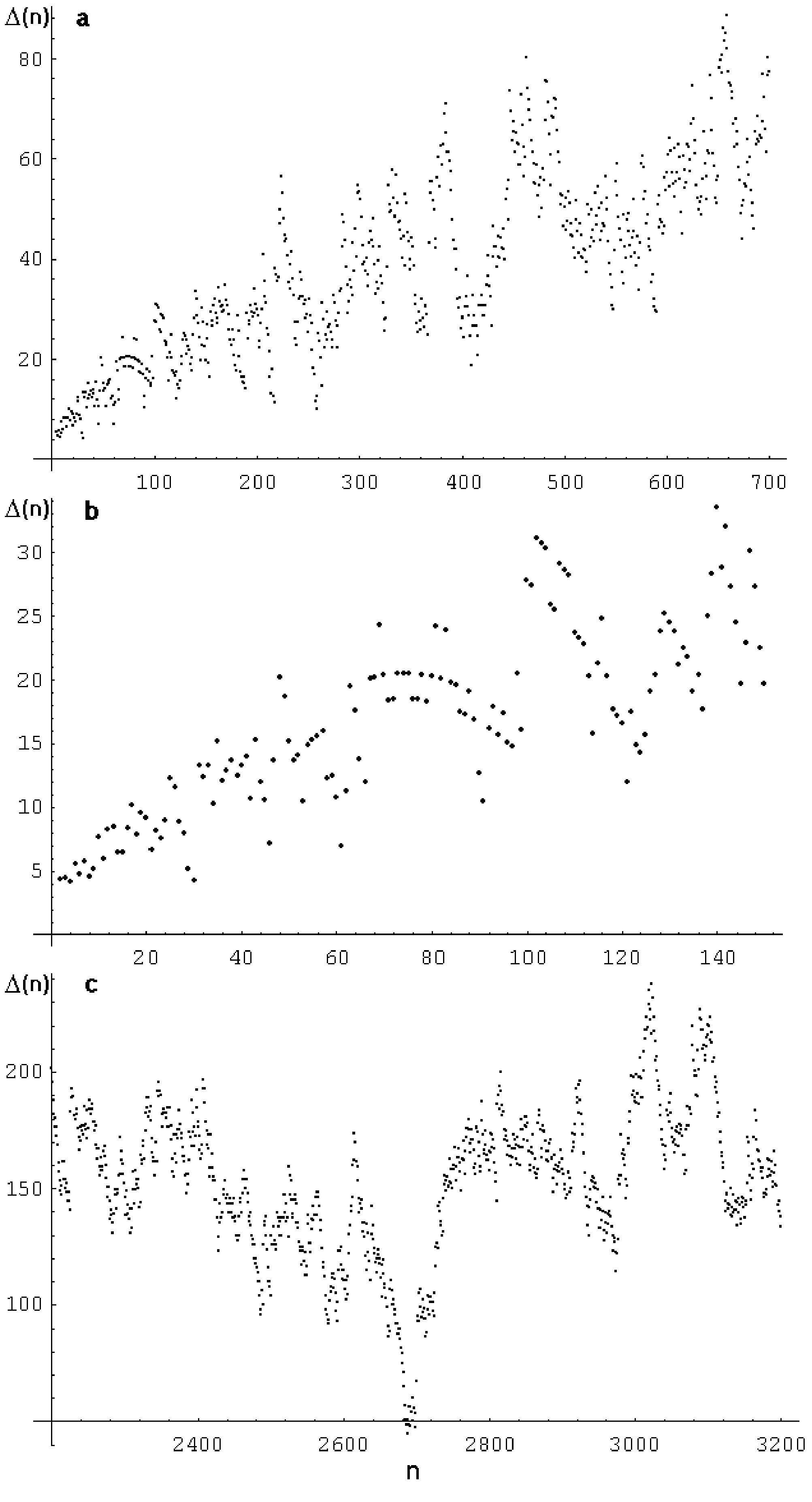

Fig. 1.

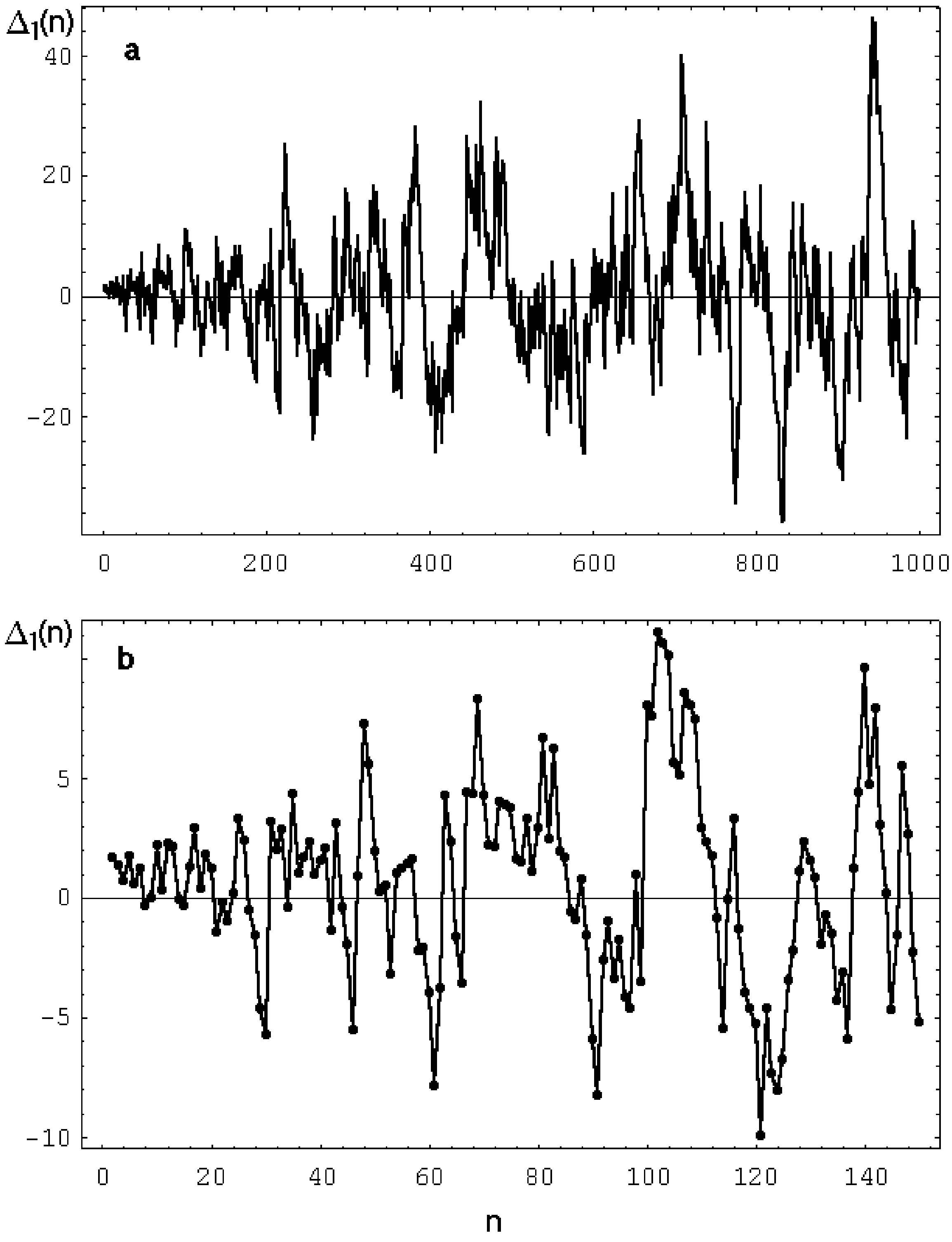

Fig. 2.

A simple analysis shows that: 1) Arcs are a nonlinear transformation of the straight lines $l_j + m_j\, n$ while subtracting correspondingly asymptotic from $p(n)$. 2) The locations of the primes on arcs or on the lengths of straight lines corresponding to them are not accidental, but reflect the presence of relations of equivalence. 3) There is a lattice for $p(n)$ (Fig. 3), which is a geometric presentation of a structure of the group of classes on reciprocally simple modulo $m_j$. Primes sequence embedded at specific sites in the lattice formed by bunches of straight lines, determined by the structure of the mentioned group of classes.

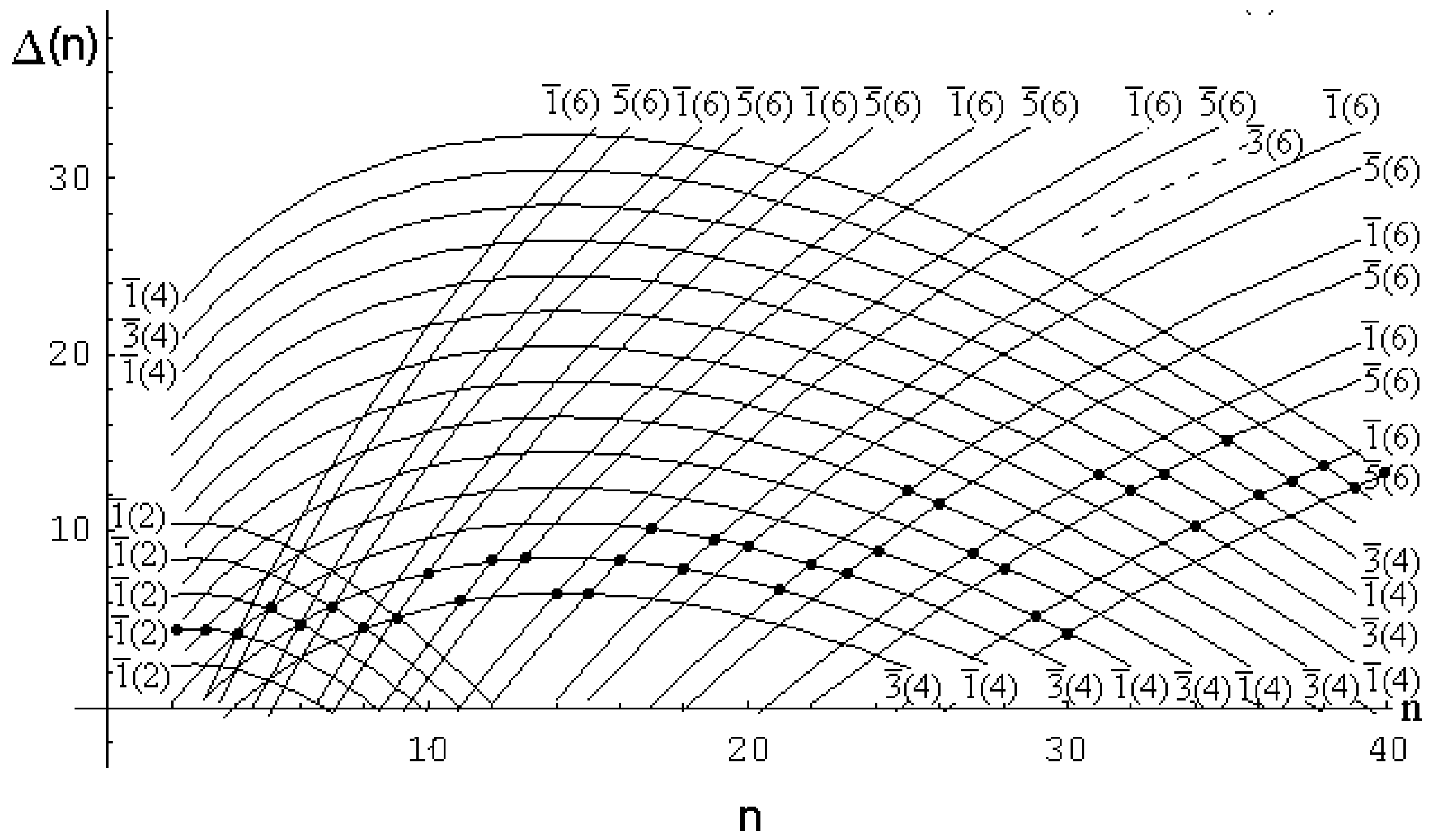

Fig. 3. A lattice $p(n)$ for begin of the sequence of natural numbers was deformed for visual demonstration. Actually it is formed by intersection of straight lines, corresponding to arithmetical progressions. The marking beside arcs means a relevant class of residues on reciprocally simple modulo. The primes lying upon the given arc belong to marked class.

So, the found arched structures of $\Delta(n)$ are described within the framework of the classical elementary numbers theory [9–12]. Their pictures can be useful in textbooks and monographs concerning problems of the distribution of primes. The visualization of the known theorems can be useful both for teaching and for researches in decryption of natural and artificial codes.

**II**. The consideration of the nature of the intermitted arches has allowed us to take a step to the genuine structure of natural sequences. For the explanation of such well organized arched structures with elements of the apparent disorder and elusive symmetry we will go back to whole sequence $N$. We will give a description of the matrix structure allowing to use it in explorations. As a matter of fact this is a potentially infinite triangular matrix with two inputs, which could be named two-dimensional sieve with the memory. One of variants of our R-matrix is given in Table 1. Here the number of symbols in $k$-th column of the triangular matrix is $k$ equal to $k$-th number of the positive integers. Each $i$-th row represents itself a periodic sequence of zeroes and units with a period $i$. Presence of unit in $i$-th row of $k$-th column points at divisibility of $k$ by $i$. Primes in this sieve are all and only those numbers $k$, which do not contain units in the corresponding column except its first and last elements, reflecting in particular, the divisibility such number $k$ only by unit and by itself. Removing a degeneration of the presentation of degrees of divisibility can be realized by increasing of a dimension of the sieve [6].

It is possible that the two-dimensional sieve with simple two letters alphabet, obviously demonstrating the structure of divisors in $N$ was considered but was not appraised at its true worth by preceding researchers including the modern number theory [13]. Periodicity along rows, periodicity along lines of divisibility, simple scheme of coincidences generating primes, twins, etc, could have been noticed three hundred years ago. Setting $a = 1$, $b = 0$ (or on the contrary), would be possible to write an equation for primes by using elementary periodic functions already well known at that time.

*Table 1*

| 1 | **2** | **3** | 4 | **5** | 6 | **7** | 8 | 9 | 10 | **11** | 12 | **13** | 14 | 15 | 16 | **17** | 18 | **19** | 20 | 21 | 22 | **23** | 24 | 25 | 26 | 27 | 28 | **29** | 30 | **31** | ... |
|---|---|---|---|---|---|---|---|---|---|---|---|---|---|---|---|---|---|---|---|---|---|---|---|---|---|---|---|---|---|---|---|
| 1 | 1 | 1 | 1 | 1 | 1 | 1 | 1 | 1 | 1 | 1 | 1 | 1 | 1 | 1 | 1 | 1 | 1 | 1 | 1 | 1 | 1 | 1 | 1 | 1 | 1 | 1 | 1 | 1 | 1 | 1 | ... |
| | 1 | **0** | 1 | **0** | 1 | **0** | 1 | 0 | 1 | **0** | 1 | **0** | 1 | 0 | 1 | **0** | 1 | **0** | 1 | 0 | 1 | **0** | 1 | 0 | 1 | 0 | 1 | **0** | 1 | **0** | ... |
| | | 1 | 0 | **0** | 1 | **0** | 0 | 1 | 0 | **0** | 1 | **0** | 0 | 1 | 0 | **0** | 1 | **0** | 0 | 1 | 0 | **0** | 1 | 0 | 0 | 1 | 0 | **0** | 1 | **0** | ... |
| | | | 1 | **0** | 0 | **0** | 1 | 0 | 0 | **0** | 1 | **0** | 0 | 0 | 1 | **0** | 0 | **0** | 1 | 0 | 0 | **0** | 1 | 0 | 0 | 0 | 1 | **0** | 0 | **0** | ... |
| | | | | 1 | 0 | **0** | 0 | 0 | 1 | **0** | 0 | **0** | 0 | 1 | 0 | **0** | 0 | **0** | 1 | 0 | 0 | **0** | 0 | 1 | 0 | 0 | 0 | **0** | 1 | **0** | ... |
| | | | | | 1 | **0** | 0 | 0 | 0 | **0** | 1 | **0** | 0 | 0 | 0 | **0** | 1 | **0** | 0 | 0 | 0 | **0** | 1 | 0 | 0 | 0 | 0 | **0** | 1 | **0** | ... |
| | | | | | | 1 | 0 | 0 | 0 | **0** | 0 | **0** | 1 | 0 | 0 | **0** | 0 | **0** | 0 | 1 | 0 | **0** | 0 | 0 | 0 | 0 | 1 | **0** | 0 | **0** | ... |
| | | | | | | | 1 | 0 | 0 | **0** | 0 | **0** | 0 | 0 | 1 | **0** | 0 | **0** | 0 | 0 | 0 | **0** | 1 | 0 | 0 | 0 | 0 | **0** | 0 | **0** | ... |
| | | | | | | | | 1 | 0 | **0** | 0 | **0** | 0 | 0 | 0 | **0** | 1 | **0** | 0 | 0 | 0 | **0** | 0 | 0 | 0 | 1 | 0 | **0** | 0 | **0** | ... |
| | | | | | | | | | 1 | **0** | 0 | **0** | 0 | 0 | 0 | **0** | 0 | **0** | 1 | 0 | 0 | **0** | 0 | 0 | 0 | 0 | 0 | **0** | 1 | **0** | ... |
| | | | | | | | | | | 1 | 0 | **0** | 0 | 0 | 0 | **0** | 0 | **0** | 0 | 0 | 1 | **0** | 0 | 0 | 0 | 0 | 0 | **0** | 0 | **0** | ... |
| | | | | | | | | | | | 1 | **0** | 0 | 0 | 0 | **0** | 0 | **0** | 0 | 0 | 0 | **0** | 1 | 0 | 0 | 0 | 0 | **0** | 0 | **0** | ... |
| | | | | | | | | | | | | 1 | 0 | 0 | 0 | **0** | 0 | **0** | 0 | 0 | 0 | **0** | 0 | 0 | 1 | 0 | 0 | **0** | 0 | **0** | ... |
| | | | | | | | | | | | | | 1 | 0 | 0 | **0** | 0 | **0** | 0 | 0 | 0 | **0** | 0 | 0 | 0 | 0 | 1 | **0** | 0 | **0** | ... |
| | | | | | | | | | | | | | | 1 | 0 | **0** | 0 | **0** | 0 | 0 | 0 | **0** | 0 | 0 | 0 | 0 | 0 | **0** | 1 | **0** | ... |
| | | | | | | | | | | | | | | | 1 | **0** | 0 | **0** | 0 | 0 | 0 | **0** | 0 | 0 | 0 | 0 | 0 | **0** | 0 | **0** | ... |
| | | | | | | | | | | | | | | | | 1 | 0 | **0** | 0 | 0 | 0 | **0** | 0 | 0 | 0 | 0 | 0 | **0** | 0 | **0** | ... |
| | | | | | | | | | | | | | | | | | 1 | **0** | 0 | 0 | 0 | **0** | 0 | 0 | 0 | 0 | 0 | **0** | 0 | **0** | ... |
| | | | | | | | | | | | | | | | | | | 1 | 0 | 0 | 0 | **0** | 0 | 0 | 0 | 0 | 0 | **0** | 0 | **0** | ... |
| | | | | | | | | | | | | | | | | | | | 1 | 0 | 0 | **0** | 0 | 0 | 0 | 0 | 0 | **0** | 0 | **0** | ... |
| | | | | | | | | | | | | | | | | | | | | 1 | 0 | **0** | 0 | 0 | 0 | 0 | 0 | **0** | 0 | **0** | ... |
| | | | | | | | | | | | | | | | | | | | | | 1 | **0** | 0 | 0 | 0 | 0 | 0 | **0** | 0 | **0** | ... |
| | | | | | | | | | | | | | | | | | | | | | | 1 | 0 | 0 | 0 | 0 | 0 | **0** | 0 | **0** | ... |
| | | | | | | | | | | | | | | | | | | | | | | | 1 | 0 | 0 | 0 | 0 | **0** | 0 | **0** | ... |
| | | | | | | | | | | | | | | | | | | | | | | | | 1 | 0 | 0 | 0 | **0** | 0 | **0** | ... |
| | | | | | | | | | | | | | | | | | | | | | | | | | 1 | 0 | 0 | **0** | 0 | **0** | ... |
| | | | | | | | | | | | | | | | | | | | | | | | | | | 1 | 0 | **0** | 0 | **0** | ... |
| | | | | | | | | | | | | | | | | | | | | | | | | | | | 1 | **0** | 0 | **0** | ... |
| | | | | | | | | | | | | | | | | | | | | | | | | | | | | 1 | 0 | **0** | ... |
| | | | | | | | | | | | | | | | | | | | | | | | | | | | | | 1 | **0** | ... |
| | | | | | | | | | | | | | | | | | | | | | | | | | | | | | | 1 | ... |
| | | | | | | | | | | | | | | | | | | | | | | | | | | | | | | | ... |

Really, let matrix element $r_{ik} = 1$ if a divisibility exactly and the remaining elements are marked by zeroes. Then it is possible to write for $r_{ik}$ (see Table 1):

$$r_{ik} = \frac{\varepsilon^2}{\varepsilon^2 + \sin(\pi k/i)}, \tag{1}$$

where ε is an infinitesimal constant. For numerical estimations at the beginning of the natural sequence it is possible to settle for value $\varepsilon^2 \sim 10^{-j}$, giving unit or zero with $j$-digit precision. The amount of nontrivial divisors in $k$-th column is

$$S_k = \sum_{i=2}^{k-1} r_{ik} , \tag{2}$$

and if $S_k = 0$, the number $k$ is a prime. Thereby, expression $S_k = 0$ is an equation for the prime number $p = k$ (see Fig. 4). Hereinafter for the number of primes less than or equal to $n$ we have:

$$\pi(n) = \sum_{k=2}^{n} \frac{\varepsilon^2}{\varepsilon^2 + S_k^2} \tag{3}$$

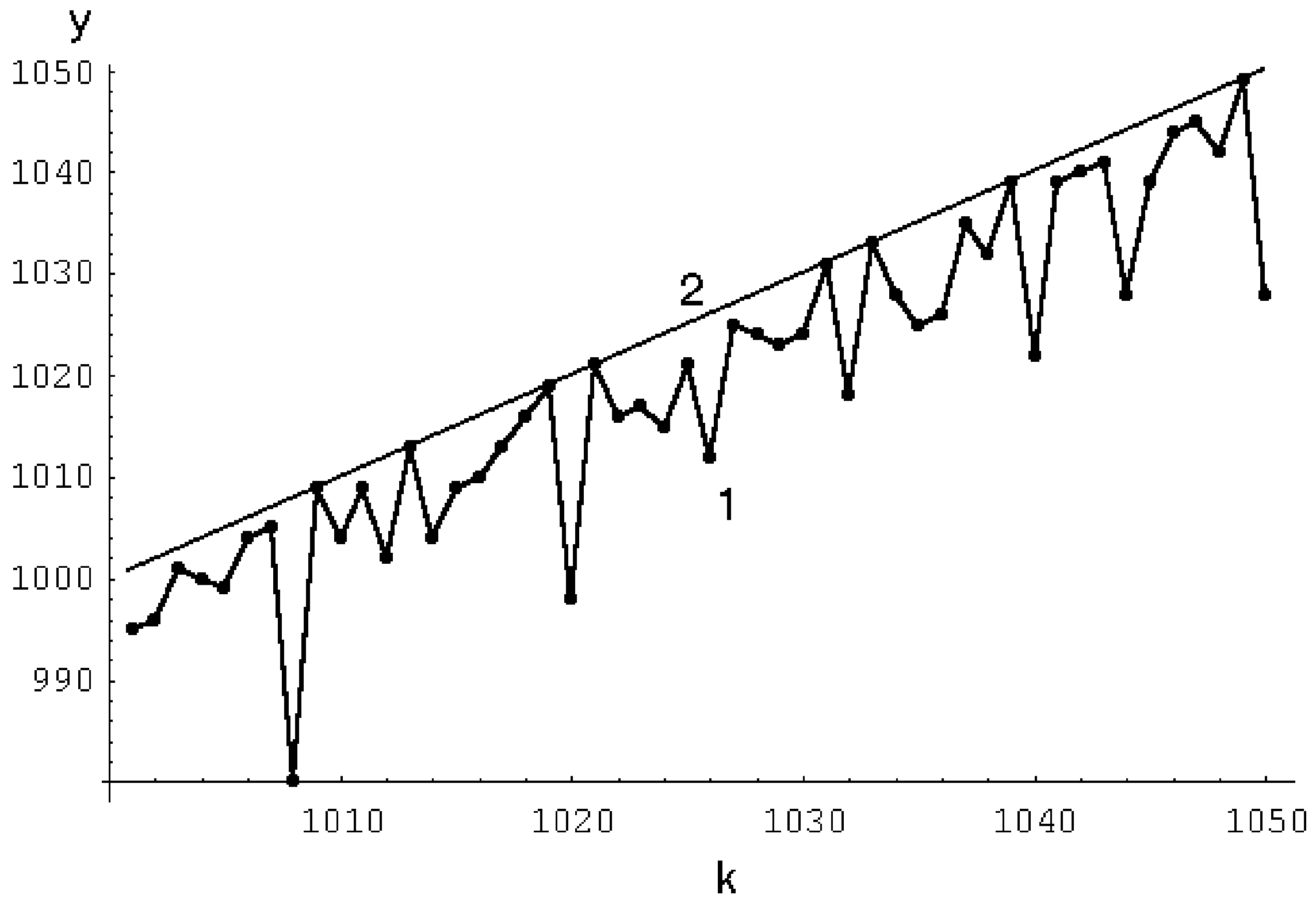

Fig. 4. Primes hanger up on the line. The function $y = k - S_k$ is shown by the curve 1, curve 2 represents $y = k$ for $k$ in the range $10001 \leq k \leq 1050$. All eight primes this range (1009, 1013, 1019, 1021, 1031, 1033, 1039, 1049) including two pair of twins, lie along a straight-line $y = k$.

Yet another expressions of $r_{ik}$ and $\pi(n)$ are:

$$r_{ik} = \lim_{\xi \to 0} \frac{1}{n!} \frac{d^n}{d\xi^n} \frac{1}{1-\xi^i} \tag{4}$$

$$\pi(n) = \sum_{k=2}^{n} \prod_{i=2}^{k-1} (1 - r_{ik}) \tag{5}$$

If need be, this formula allows to state a regular procedure of finding the average value $\pi(n)$, dispersion and others moments into interesting interval $[n]$. Obviously that if in Table 1 instead of zero to note units and accordingly instead of units to put zero, the form of the formulas a little will vary. The corresponding formulas can be recorded for twins and other consecutive primes. These formulas can be used in the search for solutions of diophantine equations. Finally, all the formulas can be represented in integral variants using the generalizing functions.

Presentation of the natural numbers sequence as a matrix formed by primitive periodic sequences in rows with peculiar superposition "generating" primes and composite numbers, allows to understand the distribution of primes that seemed up to now to be enigmatic and irregular.

It was shown [8] that the generating function for binary matrix elements in the row $m$ is $z^{m-1}/(1-z^m)$ and in the case of the color one it should be multiplied by $1/(1-z)^2$.

One can calculate the numbers of zeros $(n-1)/2 + (n-2)2/3 + \ldots$ and units $n + (n-1)/2 + (n-2)/3 \ldots$ and estimate a fraction of primes relative to $n$ for the binary matrix of the order $n$ and obtain the expression for asymptotic main term

$$\lim_{n \to \infty} \pi(n) \sim n/\mathrm{H}_n,$$

where $\mathrm{H_n}$ is the harmonic number and $\mathrm{H}_n \approx \ln n + \gamma + 1/2n + 1/12n^2 + 1/120n^4 + O(n^{-6})$. So we have $\pi(n) \sim n/\ln n$. The same result one can obtain with the use of a probability approach to estimation of event to detect the column $n$ with «coincidental rendezvois» of *(n-2)* zeros.
In reality, we must deal with only upper part (up to *n/2*) of «nilpotent» matrix.

There is a lot of works where chaotic distribution of primes in $N$ is questioned by variety of tests. These tests can not solve the problem of chaotic behavior of primes. In the end we can conclude that structure of prime distribution within $N$ is complicate, but has entirely algoritmic structure. Therefore, the above-mentioned tests are not so satisfactory.

It is necessary to note original symmetry in distribution of the prime numbers [8]. The distinctive feature of R-matrix is that anyone column (including columns meaning the prime numbers) with its elements $r_{ik}$ was preceded by the appropriate diagonal with an inclination $-45°$ with the identical set of elements (see Tabl. 1). There is also subsequent diagonal (with an inclination 45). Thus, for every prime number $p_i = n$ there is precursor in the form of diagonal going from element $[(n+1)/2; (n-1)/2]$ to element $[1; n-1]$.

**III**. The variant of a more informative sieve is given in the Table 2 where the elements of rows are numbered within each period of the appropriate sequence and the divisibility exactly is marked by zero. In this Table in an explicit form the position of primes in nodes of a lattice (Fig. 3) is reflected. For example, for $p(9)$ equal 23 and located on intersection of three arcs (of a class 1 modulus-2, a class 3 modulus-4, and a class 5 modulus-6) in Table 2 in 23-th column we discover accordingly 1 in the second row, 3 in the 4-th row and 5 in 6-th. Table 2 gives also obvious representation about multi-cyclic i.e. about a multilevel appliqué-like character of pattern of divisors in the sieve (it is seen and in Table 1). So, the pattern in first two rows repeats along $N$ with a period $P_{i=2} = 2$, for the three first rows the length of a cycle already covers $P_{i=3} = 2\times3 = 6$, and further $P_{i=4} = 2\times2\times3 = 12$, $P_{i=5} = 2\times2\times3\times5 = 60$, $P_{i=6} = 60$, $P_{i=7} = 420$, $P_{i=8} = 840$, ..., $P_{i=41} = 219060189739591200$, $P_{i=43} = 9419588158802421600$, …. On the average length of a cycle $P_i$ grows with its «depth» $i$ as Exp($i$), but always the cycle of a higher rank contains an integer of cycles of all lower ranks. We should like to note a curious fact: $P_{i=43} + p(21)$, $P_{i=43} + p(22)$, $P_{i=43} + p(23)$, $P_{i=43} + p(24)$ are too primes.

Turning back to the physical aspect of the problem, we will indicate orders of values for great systems, for which the evaluation of the behavior $\pi(n)$ and $p(n)$ as a whole is interesting. In the most ambitious physical theories of the type "great unification" relations of the age of Universe to Planck time-scale and scale of observable Universe to the Planck length $t_{un}/t_{pl} = l_{un}/l_{pl} \sim 10^{61}$ play a fundamental role. One can obtain the same order of magnitude for the relations of the mass Universe (in models with dark matter [11]) to greatly possible mass of elementary particles. The value $1 - 1/10^{61}$ is of the same interest for us. The value of the logarithm is nearly twice less for «terrestrial touchstone» experiments in order to feel out «fuzzy» or «foamy» structure of the space-time [12]. For more modest physical and physicist-chemical problems an order of values is even less, unless speak about the combinatory evaluations.

*Table 2*

| 1 | **2** | **3** | 4 | **5** | 6 | **7** | 8 | 9 | 10 | **11** | 12 | **13** | 14 | 15 | 16 | **17** | 18 | **19** | 20 | 21 | 22 | **23** | 24 | 25 | 26 | 27 | 28 | **29** | 30 | **31** | ... |
|---|---|---|---|---|---|---|---|---|---|---|---|---|---|---|---|---|---|---|---|---|---|---|---|---|---|---|---|---|---|---|---|
| **0** | **0** | **0** | 0 | **0** | 0 | **0** | 0 | 0 | 0 | **0** | 0 | **0** | 0 | 0 | 0 | **0** | 0 | **0** | 0 | 0 | 0 | **0** | 0 | 0 | 0 | 0 | 0 | **0** | 0 | **0** | ... |
| 1 | **0** | **1** | 0 | **1** | 0 | **1** | 0 | 1 | 0 | **1** | 0 | **1** | 0 | 1 | 0 | **1** | 0 | **1** | 0 | 1 | 0 | **1** | 0 | 1 | 0 | 1 | 0 | **1** | 0 | **1** | ... |
| 1 | 2 | **0** | 1 | **2** | 0 | **1** | 2 | 0 | 1 | **2** | 0 | **1** | 2 | 0 | 1 | **2** | 0 | **1** | 2 | 0 | 1 | **2** | 0 | 1 | 2 | 0 | 1 | **2** | 0 | **1** | ... |
| 1 | 2 | 3 | 0 | **1** | 2 | **3** | 0 | 1 | 2 | **3** | 0 | **1** | 2 | 3 | 0 | **1** | 2 | **3** | 0 | 1 | 2 | **3** | 0 | 1 | 2 | 3 | 0 | **1** | 2 | **3** | ... |
| 1 | 2 | 3 | 4 | **0** | 1 | **2** | 3 | 4 | 0 | **1** | 2 | **3** | 4 | 0 | 1 | **2** | 3 | **4** | 0 | 1 | 2 | **3** | 4 | 0 | 1 | 2 | 3 | **4** | 0 | **1** | ... |
| 1 | 2 | 3 | 4 | 5 | 0 | **1** | 2 | 3 | 4 | **5** | 0 | **1** | 2 | 3 | 4 | **5** | 0 | **1** | 2 | 3 | 4 | **5** | 0 | 1 | 2 | 3 | 4 | **5** | 0 | **1** | ... |
| 1 | 2 | 3 | 4 | 5 | 6 | **0** | 1 | 2 | 3 | **4** | 5 | **6** | 0 | 1 | 2 | **3** | 4 | **5** | 6 | 0 | 1 | **2** | 3 | 4 | 5 | 6 | 0 | **1** | 2 | **3** | ... |
| 1 | 2 | 3 | 4 | 5 | 6 | 7 | 0 | 1 | 2 | **3** | 4 | **5** | 6 | 7 | 0 | **1** | 2 | **3** | 4 | 5 | 6 | **7** | 0 | 1 | 2 | 3 | 4 | **5** | 6 | **7** | ... |
| 1 | 2 | 3 | 4 | 5 | 6 | 7 | 8 | 0 | 1 | **2** | 3 | **4** | 5 | 6 | 7 | **8** | 0 | **1** | 2 | 3 | 4 | **5** | 6 | 7 | 8 | 0 | 1 | **2** | 3 | **4** | ... |
| 1 | 2 | 3 | 4 | 5 | 6 | 7 | 8 | 9 | 0 | **1** | 2 | **3** | 4 | 5 | 6 | **7** | 8 | **9** | 0 | 1 | 2 | **3** | 4 | 5 | 6 | 7 | 8 | **9** | 0 | **1** | ... |
| 1 | 2 | 3 | 4 | 5 | 6 | 7 | 8 | 9 | 10 | **0** | 1 | **2** | 3 | 4 | 5 | **6** | 7 | **8** | 9 | 10 | 0 | **1** | 2 | 3 | 4 | 5 | 6 | **7** | 8 | **9** | ... |
| 1 | 2 | 3 | 4 | 5 | 6 | 7 | 8 | 9 | 10 | 11 | 0 | **1** | 2 | 3 | 4 | **5** | 6 | **7** | 8 | 9 | 10 | **11** | 0 | 1 | 2 | 3 | 4 | **5** | 6 | **7** | ... |
| 1 | 2 | 3 | 4 | 5 | 6 | 7 | 8 | 9 | 10 | 11 | 12 | **0** | 1 | 2 | 3 | **4** | 5 | **6** | 7 | 8 | 9 | **10** | 11 | 12 | 0 | 1 | 2 | **3** | 4 | **5** | ... |
| 1 | 2 | 3 | 4 | 5 | 6 | 7 | 8 | 9 | 10 | 11 | 12 | 13 | 0 | 1 | 2 | **3** | 4 | **5** | 6 | 7 | 8 | **9** | 10 | 11 | 12 | 13 | 0 | **1** | 2 | **3** | ... |
| 1 | 2 | 3 | 4 | 5 | 6 | 7 | 8 | 9 | 10 | 11 | 12 | 13 | 14 | 0 | 1 | **2** | 3 | **4** | 5 | 6 | 7 | **8** | 9 | 10 | 11 | 12 | 13 | **14** | 0 | **1** | ... |
| 1 | 2 | 3 | 4 | 5 | 6 | 7 | 8 | 9 | 10 | 11 | 12 | 13 | 14 | 15 | 0 | **1** | 2 | **3** | 4 | 5 | 6 | **7** | 8 | 9 | 10 | 11 | 12 | **13** | 14 | **15** | ... |
| 1 | 2 | 3 | 4 | 5 | 6 | 7 | 8 | 9 | 10 | 11 | 12 | 13 | 14 | 15 | 16 | **0** | 1 | **2** | 3 | 4 | 5 | **6** | 7 | 8 | 9 | 10 | 11 | **12** | 13 | **14** | ... |
| 1 | 2 | 3 | 4 | 5 | 6 | 7 | 8 | 9 | 10 | 11 | 12 | 13 | 14 | 15 | 16 | 17 | 0 | **1** | 2 | 3 | 4 | **5** | 6 | 7 | 8 | 9 | 10 | **11** | 12 | **13** | ... |
| 1 | 2 | 3 | 4 | 5 | 6 | 7 | 8 | 9 | 10 | 11 | 12 | 13 | 14 | 15 | 16 | 17 | 18 | **0** | 1 | 2 | 3 | **4** | 5 | 6 | 7 | 8 | 9 | **10** | 11 | **12** | ... |
| 1 | 2 | 3 | 4 | 5 | 6 | 7 | 8 | 9 | 10 | 11 | 12 | 13 | 14 | 15 | 16 | 17 | 18 | 19 | 0 | 1 | 2 | **3** | 4 | 5 | 6 | 7 | 8 | **9** | 10 | **11** | ... |
| 1 | 2 | 3 | 4 | 5 | 6 | 7 | 8 | 9 | 10 | 11 | 12 | 13 | 14 | 15 | 16 | 17 | 18 | 19 | 20 | 0 | 1 | **2** | 3 | 4 | 5 | 6 | 7 | **8** | 9 | **10** | ... |
| 1 | 2 | 3 | 4 | 5 | 6 | 7 | 8 | 9 | 10 | 11 | 12 | 13 | 14 | 15 | 16 | 17 | 18 | 19 | 20 | 21 | 0 | **1** | 2 | 3 | 4 | 5 | 6 | **7** | 8 | **9** | ... |
| 1 | 2 | 3 | 4 | 5 | 6 | 7 | 8 | 9 | 10 | 11 | 12 | 13 | 14 | 15 | 16 | 17 | 18 | 19 | 20 | 21 | 22 | **0** | 1 | 2 | 3 | 4 | 5 | **6** | 7 | **8** | ... |
| 1 | 2 | 3 | 4 | 5 | 6 | 7 | 8 | 9 | 10 | 11 | 12 | 13 | 14 | 15 | 16 | 17 | 18 | 19 | 20 | 21 | 22 | 23 | 0 | 1 | 2 | 3 | 4 | **5** | 6 | **7** | ... |
| 1 | 2 | 3 | 4 | 5 | 6 | 7 | 8 | 9 | 10 | 11 | 12 | 13 | 14 | 15 | 16 | 17 | 18 | 19 | 20 | 21 | 22 | 23 | 24 | 0 | 1 | 2 | 3 | **4** | 5 | **6** | ... |
| 1 | 2 | 3 | 4 | 5 | 6 | 7 | 8 | 9 | 10 | 11 | 12 | 13 | 14 | 15 | 16 | 17 | 18 | 19 | 20 | 21 | 22 | 23 | 24 | 25 | 0 | 1 | 2 | **3** | 4 | **5** | ... |
| 1 | 2 | 3 | 4 | 5 | 6 | 7 | 8 | 9 | 10 | 11 | 12 | 13 | 14 | 15 | 16 | 17 | 18 | 19 | 20 | 21 | 22 | 23 | 24 | 25 | 26 | 0 | 1 | **2** | 3 | **4** | ... |
| 1 | 2 | 3 | 4 | 5 | 6 | 7 | 8 | 9 | 10 | 11 | 12 | 13 | 14 | 15 | 16 | 17 | 18 | 19 | 20 | 21 | 22 | 23 | 24 | 25 | 26 | 27 | 0 | **1** | 2 | **3** | ... |
| 1 | 2 | 3 | 4 | 5 | 6 | 7 | 8 | 9 | 10 | 11 | 12 | 13 | 14 | 15 | 16 | 17 | 18 | 19 | 20 | 21 | 22 | 23 | 24 | 25 | 26 | 27 | 28 | **0** | 1 | **2** | ... |
| 1 | 2 | 3 | 4 | 5 | 6 | 7 | 8 | 9 | 10 | 11 | 12 | 13 | 14 | 15 | 16 | 17 | 18 | 19 | 20 | 21 | 22 | 23 | 24 | 25 | 26 | 27 | 28 | 29 | 0 | **1** | ... |
| 1 | 2 | 3 | 4 | 5 | 6 | 7 | 8 | 9 | 10 | 11 | 12 | 13 | 14 | 15 | 16 | 17 | 18 | 19 | 20 | 21 | 22 | 23 | 24 | 25 | 26 | 27 | 28 | 29 | 30 | **0** | ... |
| | | | | | | | | | | | | | | | | | | | | | | | | | | | | | | | … |

**Acknowledgements.** Authors are grateful to A. G. Bashkirov and V. V. Guchuk for discussing the results of investigations and to the Management of the Institute for Dynamics of Geospheres RAS for giving possibility to carry it out.